\documentclass[12pt]{article}
\usepackage{amsmath,amsthm,amsfonts,amssymb,amscd}
\usepackage[usenames]{color}
\usepackage[dvips]{graphicx}
\usepackage{epsfig}
\def\phi{\varphi}

\begin{document}

\newtheorem{lem}{Lemma}
\newtheorem{predl}{Proposition}
\newtheorem{theorem}{Theorem}
\newtheorem{defin}{Definition}
\newtheorem{zam}{Remark}

\bigskip

\centerline{\bf TESTS OF EXPONENTIALITY BASED ON ARNOLD-VILLASENOR}
\centerline{\bf CHARACTERIZATION, AND THEIR EFFICIENCIES}

\bigskip

\bigskip

\centerline{ \it  M. Jovanovi\'c $^{a}$, B. Milo\v sevi\'c $^{a}$, Ya. Yu. Nikitin $^{b,}$\footnote{Research of Ya. Yu. Nikitin and K. Yu. Volkova was supported by RFBR grant No. 13-01-00172, and by SPbGU grant No. 6.38.672.2013},}
\centerline{\it M. Obradovi\'c $^{a}$, K. Yu. Volkova $^{c}$}

\bigskip
\medskip

$^a$Faculty of Mathematics, University of Belgrade, Studenski trg 16, Belgrade, Serbia;

\medskip

$^b$Department of Mathematics and Mechanics, Saint-Petersburg State University, Universitetsky pr. 28,
Stary Peterhof 198504, Russia, and National Research University - Higher School of Economics, Souza
Pechatnikov, 16, St.Peters\-burg 190008, Russia;

\medskip

$^c$ Department of Mathematics and Mechanics, Saint-Petersburg State University, Universitetsky pr. 28,
Stary Peterhof 198504, Russia.

\bigskip

\noindent {\bf Abstract.}{ \small We propose two families of scale-free exponentiality tests based on the recent characterization of exponentiality by Arnold and Villasenor. The test statistics are based on suitable functionals of $U$-empirical distribution functions. The family of integral statistics can be reduced to $V$- or $U$-statistics with relatively simple non-degenerate kernels. They are asymptotically normal and have reasonably high local Bahadur efficiency under common alternatives. This efficiency is compared with simulated powers of new tests.
On the other hand, the Kolmogorov type tests demonstrate very low local Bahadur efficiency and rather moderate power for common alternatives, and can hardly be recommended to practitioners. We also explore the conditions of local asymptotic optimality of new tests and describe for both families special "most favorable" alternatives for which the  tests are fully efficient.}

\bigskip


\bigskip


 \section{ Introduction}

Exponential distribution plays an essential role in Probability and Statistics
since various models with exponentially distributed observations often appear
in applications such as survival analysis, reliability theory, engineering,
demography, etc. Therefore, testing exponentiality is one of the most important problems
in goodness-of-fit theory.

There exists a multitude of tests for this problem which are based on various
ideas (see books and reviews
 \cite{Ahsan}, \cite{A90}, \cite{BB95}, \cite{DY}, \cite{HM02},
 \cite{NJ02}). Among them many tests are based on characterizations
of exponential law, in particular on loss-of-memory property (\cite{AA},
\cite{Ang}, \cite{Kou}, \cite{Kou2}, \cite{Niknik}) and some other characterizations (\cite{BarHen}, \cite{JansV}, \cite{Litv2}, \cite{NikVol}, \cite{Nou}, \cite{Rank}, \cite{RT}). The construction of tests based on characterizations is a relatively fresh idea which gradually becomes one of main directions in goodness-of-fit testing.

In this paper we present new tests for exponentiality based on Arnold-Villasenor characterization.
In \cite{Ar} Arnold and Villasenor stated the following hypothesis:

\emph{Let $\cal F$ be the class of distributions whose densities have derivatives of all orders in the neighbourhood of zero and let $X_1, X_{2}, \ldots, X_n$ be non-negative independent identically
 distributed (i.i.d.) random variables with distribution function (d.f.) $F$ from class $\cal F.$
 Then the random variables  $\max(X_1, X_2, \ldots, X_k)$ and $\sum_{i=1}^{k}\frac{X_i}{i}$
are equally distributed if and only if the d.f. $F$ is exponential.}

They were able to prove this hypothesis only for $k=2$. Later Yanev and Chakraborty in \cite{Yanev} proved that this hypothesis was also true for $k=3$. We think that the validity of Arnold-Villasenor hypothesis is very likely,
and it will be proved in the nearest future. This is sustained by the fact that recently Chakraborty and Yanev proved the correctness of the related hypothesis from \cite{Ar} for any $k$ (see details in \cite{Chakra}).

Let $X_1,X_2,\ldots,X_n$ be i.i.d.
observations having the continuous d.f. $F$ from the class $\cal F$. We are testing the composite hypothesis
of exponentiality {\it $H_{0}: F(x)$ belongs to exponential family of distributions $ {\cal {E(\lambda)}}$  with the density $
f(x)=\lambda e^{-\lambda x}, x \geq 0,$} where $\lambda>0$ is an unknown parameter.

Let $F_n(t)=n^{-1}\sum_{i=1}^n\textbf{1}\{X_i<t\}, t \in \mathbb{R},$ be the usual empirical d.f.
based on the observations $X_1,X_{2},\ldots,X_n.$ In compliance with Arnold-Villasenor characterization for $t\geq 0$ we introduce
 the so-called $V$-empirical d.f.'s (see \cite{Jan}, \cite{Kor}) according to the formulae
\begin{eqnarray*}
H_n^{(k)}(t)&=&\frac{1}{n^{k}}\sum_{i_1,i_{2}, \ldots,i_k
=1}^{n}\textbf{1}\{\max(X_{i_1}, X_{i_2}, \ldots, X_{i_k})<t\},\\
G_n^{(k)}(t)&=&\frac{1}{n^{k}k!}\sum_{i_1, \ldots,i_k =1}^{n}
\bigg[\sum_{\pi(j_1, \ldots,j_k)}
\textbf{1}\{\frac{X_{i_1}}{j_1}+\frac{X_{i_2}}{j_2}+\ldots+\frac{X_{i_k}}{j_k}<
t\}\bigg],
\end{eqnarray*}
\noindent where  $\pi(j_1, \ldots,j_k)$ represents the set of
all $k!$ permutations of natural numbers $1, 2, \ldots, k,$ $k\geq 2$.

It is well-known that the properties of $V$- and $U$-empirical d.f.'s are similar to
the properties of usual empirical d.f.'s. In particular, Glivenko-Cantelli theorem is valid in this case (see \cite{HJS}, \cite{Jan}). Hence, according to Arnold-Villasenor characterization, the empirical d.f.'s  $H_n^{(k)}$ and  $G_n^{(k)}$ should be close for large $n$  under $H_0$, and we can
measure their proximity using appropriate test statistics.

Let us introduce two new sequences of statistics depending on natural $k>1$ which are invariant with respect to
the scale parameter $\lambda:$
\begin{align}
I_n^{(k)}&=\int_{0}^{\infty} \left(H_n^{(k)}(t)-G_n^{(k)}(t)\right)dF_n(t),\label{I_n_Ar}\\
D_n^{(k)}&=\sup_{t \geq 0}\mid H_n^{(k)}(t)-G_n^{(k)}(t) \mid\label{D_n_Ar},
\end{align}
\noindent where $k\geq 2$.

Large values of $I_n^{(k)}$ and $D_n^{(k)}$ are significant for rejection of null hypothesis. The sequence of statistics $I_n^{(k)}$ is not always consistent but nevertheless the  consistency takes place for many common alternatives. At first glance the sequence of statistics of omega-square type
$$
W_n^{(k)} = \int_{0}^{\infty} \left(H_n^{(k)}(t)-G_n^{(k)}(t)\right)^2 dF_n(t),
$$
could seem more adequate choice, but their asymptotic theory is very intricate and is currently underdeveloped. In the same time the statistics $I_n^{(k)}$ are usually asymptotically normal. As to the sequence $D_n^{(k)},$  it is consistent for any alternative.

 In what follows we describe the limiting distributions and large deviations
 of both sequences of statistics under  $H_0$, and calculate their
 local Bahadur efficiency under different alternatives. We also
 analyze the conditions of local asymptotic optimality of new statistics.
 In this regard we refer to the results from the theory of
 $U$- and $V$-statistics and the theory of Bahadur efficiency
 (\cite{Bahadur},  \cite{anirban}, \cite{Kor}, \cite{Nik}).

We have selected the Bahadur approach as a method of calculation
of asymptotic efficiency for our tests because the Kolmogorov-type
statistics $D_n^{(k)}$ are not asymptotically normal under null-hypothesis,
and therefore the classical Pitman approach is not applicable.
In case of integral statistic $I_n^{(k)},$ local Bahadur efficiency
and Pitman efficiency coincide (\cite{Bah60}, \cite{Wie}).

We supplement our research with simulated powers
 which principally  support  the theoretical values of efficiency.

\section {Integral statistic  $I_{n}^{(k)}$}

Without loss of generality we can assume that  $\lambda=1$. The statistic
$I_{n}^{(k)}$ is asymptotically equivalent to the $V$-statistic of degree $(k+1)$ with
the centered kernel $\Psi_k (X_1, X_2,\ldots, X_{k+1})$ given by
\begin{eqnarray*}
\Psi_k (X_1, X_2,\ldots, X_{k+1})&=&\frac{1}{k+1}\bigg[\sum\limits_{i=1}^{k+1}\textbf{1}\{\max(X_{1}, \ldots,X_{i-1}, X_{i+1}\ldots, X_{k+1})<X_{i}\} \\
&-&\frac{1}{k!}\sum\limits_{i=1}^{k+1}\! \sum_{\pi(j_1,\ldots,j_k)}\!\!\! \textbf{1}\{\frac{X_{1}}{j_1}+\!\ldots\!+\frac{X_{i-1}}{j_{i-1}}+\frac{X_{i+1}}{j_{i+1}}+\!\ldots\!+\frac{X_{k+1}}{j_{k+1}}<
X_{i}\}\bigg].
\end{eqnarray*}

It is well-known that non-degenerate $U$-
and $V$-statistics are asymptotically normal (\cite{Hoeffding},
\cite{Kor}). To show that the kernel $\Psi_k (X_1, X_2,\ldots,
X_{k+1})$ is non-degenerate, let us calculate its projection $\psi_k (s)$ under null hypothesis. For
fixed $X_{k+1}=s$ this projection has the form:
\begin{eqnarray*}
\psi_k (s)&=&E(\Psi_k (X_1, X_2,\ldots, X_{k+1})
|X_{k+1}=s)=\frac{1}{k+1}P(\max(X_1,\ldots,X_k)<s)\\
&+&\frac{k}{k+1}P(\max(s,X_2,\ldots,X_k)<X_1)
-\frac{1}{(k+1)!} \sum_{\pi(j_1, \ldots,j_k)}
P(\frac{X_1}{j_1}+\ldots+\frac{X_k}{j_k}< s)\\
&-& \frac{k}{(k+1)!}
\sum_{\pi(j_1, \ldots,j_k)}P(\frac{s}{j_1}+\frac{X_2}{j_2}+\ldots+\frac{X_k}{j_k}<
X_1).
\end{eqnarray*}

\noindent It follows from Arnold and Villasenor's characterization that the first and the third term in the right hand side coincide, so they cancel out.

 Next we calculate the second term:
\begin{eqnarray*}
\frac{k}{k+1}P(\max(s,X_2,\ldots,X_k)<X_1)&=&\frac{k}{k+1}\int_0^{\infty} \textbf{1}\{s<t\}P(X_2<t,
 \ldots, X_k<t)dF(t)\\
 &=&\frac{k}{k+1}\int_s^{\infty}F^{k-1}(s)dF(s)= \frac{1}{k+1}\left(1-F^k (s)\right),
\end{eqnarray*}
\noindent where $F(x)=1-e^{-x}$. It remains to calculate the last term. Since
\begin{eqnarray*}
P(\frac{s}{j_1} +\frac{X_2}{j_2}+\ldots+\frac{X_k}{j_k}<
X_1) &=& \int_{0}^{\infty}e^{-x_2}dx_2\ldots\int_{0}^{\infty}e^{-x_k}dx_k
\int\limits_{\frac{s}{j_1}+\frac{x_2}{j_2}+\ldots+\frac{x_k}{j_k}}^{\infty}e^{-x_1}dx_1\\
 &=& \frac{1}{(k+1)}(1+\frac{1}{j_1})e^{-s/{j_1}},
\end{eqnarray*}
after summing this expression over all permutations of indices $j_1,j_{2}, \ldots,j_k$ and some additional calculations, we get that the fourth term is
$\frac{1}{(k+1)^2}\sum_{r=1}^{k}(1+\frac{1}{r})e^{-s/{r}}.$

Finally we obtain the following expression for the projection $\psi_k$ of the kernel $\Psi_k :$
\begin{gather}
\psi_k (s)=
\frac{1-(1- e^{-s})^k}{k+1}-\frac{1}{(k+1)^2}\sum_{r=1}^{k}(1+\frac{1}{r})e^{-s/{r}}.\label{psi_Ar}
\end{gather}

\noindent It is easy to show that $E(\psi_k (X_1))=0$. After some calculations we get that
the variance of this projection is
\begin{eqnarray}\label{sigma}
\nonumber \Delta_k^2 &=&{\rm Var}(\psi_k (X_1))=\int\limits_{0}^\infty \psi_k^2(s)e^{-s}ds=\frac{1}{(k+1)^{3}}\bigg[\frac{-12k^{4}-38k^{3}-35k^{2}-11k}{4(k+1)^{2}(k+2)(2k+1)}\\
&+&2k!\sum\limits_{r=1}^{k}\frac{1}{(k+1+\frac{1}{r})(k+\frac{1}{r})
\cdots(2+\frac{1}{r})}+\frac{2}{k+1}\sum\limits_{1\leq i<j\leq k}\frac{1}{i+j+ij}\bigg].
\end{eqnarray}

\noindent It is clear from (\ref{sigma}) that the kernel $\Psi_k$ is non-degenerate for any $k.$

 In fact if the kernel is non-degenerate, we can consider instead of $V$-statistic $I_n^{(k)}$ the corresponding  $U$-statistic  with the same kernel which has  very similar asymptotic properties but is considerably simpler for calculation.

\subsection{Local Bahadur efficiency}

Let $G(\cdot,\theta)$, $\theta \geq 0$, be a family of d.f.'s with densities $g(\cdot,\theta)$, such that $G(\cdot,0)\in \mathcal{E}(\lambda)$.
The measure of Bahadur efficiency (BE) for any sequence $\{T_n\}$ of test statistics is the exact slope
$c_{T}(\theta)$ describing the rate of exponential decrease for the
attained level under the alternative d.f. $G(\cdot,\theta)$, $\theta > 0$. According to Bahadur theory  (\cite{Bahadur}, \cite{Nik}) the exact slopes may be found by
using the following proposition.

\noindent {\bf Proposition }\,{\it Suppose that the following two
conditions hold:
\[
\hspace*{-3.5cm} \mbox{a)}\qquad  T_n \
\stackrel{\mbox{\scriptsize $P_\theta$}}{\longrightarrow} \
b(\theta),\qquad \theta > 0,\nonumber \] where $-\infty <
b(\theta) < \infty$, and $\stackrel{\mbox{\scriptsize
$P_\theta$}}{\longrightarrow}$ denotes convergence in probability
under $G(\cdot, \theta)$.
\[
\hspace*{-2cm} \mbox{b)} \qquad \lim_{n\to\infty} n^{-1} \ln
P_{H_0} \left( T_n \ge t  \right)  =  - h(t)  \nonumber
\]
for any $t$ in an open interval $I,$ on which $h$ is
continuous and $\{b(\theta), \: \theta > 0\}\subset I$. Then
$c_T(\theta) = 2h(b(\theta)).$}

The exact slopes always satisfy the inequality (\cite{Bahadur}, \cite{Nik})
\begin{equation}
\label{Ragav}
c_T(\theta) \leq 2 K(\theta),\, \theta > 0,
\end{equation}
where $K(\theta)$ is the Kullback-Leibler "distance" between the alternative $H_1$ and the null hypothesis $H_0.$ In our case $H_0$ is composite, hence for any
alternative density $g(x,\theta)$ one has
\begin{equation}\label{kul}
K(\theta) = \inf_{\lambda>0} \int_0^{\infty} \ln [g(x,\theta) / \lambda \exp(-\lambda x) ] g(x,\theta) \ dx.
\end{equation}
This quantity can be easily calculated as $\theta \to 0$ for particular alternatives.
According to (\ref{Ragav}), the local BE of the sequence of statistics ${T_n}$ is defined as
$$
e^B (T) = \lim_{\theta \to 0} \frac{c_T(\theta)}{2K(\theta)}.
$$

\subsection{Integral statistic $I_{n}^{(2)}$}

For $k=2$ from (\ref{psi_Ar}) and (\ref{sigma}) we get that the projection of the kernel $\Psi_{2}(X,Y,Z)$ is equal to
\begin{gather}\label{psi_1}
\psi_2 (s)=\frac{4}{9}e^{-s}-\frac{1}{3}e^{-2s}-\frac{1}{6}e^{-s/2}
\end{gather}
\noindent and its variance is
\begin{gather*}
\Delta_2^2 = \int_{0}^{\infty} \psi_2^2 (s) e^{-s}ds =\frac{5}{13608}\approx
0.000367.
\end{gather*}

 Applying Hoeffding's theorem for U-statistics with non-degenerate kernels (see \cite{Hoeffding}, \cite{Kor})), as $n \rightarrow \infty,$ we obtain
$$\sqrt{n}I_n^{(2)} \stackrel{d}{\longrightarrow}{\cal
{N}}\left(0,\frac{5}{1512}\right).$$

Let us now find the logarithmic asymptotics of large deviations of the sequence
of statistics $I_{n}^{(2)}$ under null hypothesis. The kernel $\Psi_2 $ is centered, non-degenerate and bounded.
Applying the results on large deviations of non-degenerate $U$- and $V$-statistics from \cite{nikiponi} (see also \cite{anirban}, \cite{Niki10}),  we state the following theorem:
\begin{theorem}
 For $a>0$ it holds
$$
\lim_{n\to \infty} n^{-1} \ln P_{H_0} ( I_n^{(2)}  >a) = - f (a),
$$
where the function $f$ is analytic for sufficiently small $a>0,$ moreover
\begin{gather}
f(a) \sim \frac{a^2}{18 \Delta_2^2}=  \frac{756}{5}a^2=151.2 a^2, \,
\, \mbox{as} \, \, a \to 0.\label{th1}
\end{gather}
\end{theorem}

According to the law of large numbers
for $U$- and $V$-statistics (\cite{Kor}), the limit in probability under alternative $H_1$ is equal to
\begin{gather*}
b_I^{(2)}(\theta)=P_{\theta}(\max(X,Y)<Z)-P_{\theta}(X+\frac{Y}{2}< Z).
\end{gather*}

It is easy to show (see also \cite{NiPe}), that
\begin{equation}\label{bI2}
b_I^{(2)}(\theta) \sim 3\theta  \int_{0}^{\infty} \psi_2 (s)h(s)ds,\;\;\; \theta \rightarrow 0,
\end{equation}
where $h(x)=\frac{\partial}{\partial\theta}g_1(x,\theta)\mid _{\theta=0}$ and $\psi_2(s)$ is the projection from (\ref{psi_1}).

We present the following common alternatives against exponentiality which will be considered for all tests in this paper:

\begin{enumerate}

\item[i)] Makeham distribution with the density
$$g_1(x,\theta)=(1+\theta(1-e^{-x}))\exp(-x-\theta( e^{-x}-1+x)),\theta > 0, x\geq 0;$$

\item[ii)] Weibull distribution with the density
$$g_2(x,\theta)=(1+\theta)x^\theta \exp(-x^{1+\theta}),\theta > 0, x\geq 0;$$

\item[iii)] gamma distribution with the density
$$g_3(x,\theta)=\frac{x^{\theta}}{\Gamma(\theta+1)}e^{-x}, \theta > 0, x\geq 0;$$


\item[iv)] exponential mixture with negative weights (EMNW($\beta$)) (see \cite{vjevremovic})
$$ g_4(x)=(1+\theta)e^{-x}-\theta\beta e^{-\beta x}, x\geq 0, \theta\in \big(0,\frac{1}{\beta-1}\big]$$
\end{enumerate}

Let us calculate the local Bahadur efficiencies for these alternatives.

For the Makeham alternative from \eqref{bI2} we get that
\begin{eqnarray*}
b_I^{(2)}(\theta) &\sim& 3\theta
\int_{0}^{\infty}(\frac{4}{9}e^{-s}-\frac{1}{3}e^{-2s}-\frac{1}{6}e^{-s/2})
e^{-s}(2-2e^{-s}-s)ds\\
 &=&\frac{\theta}{90} \approx
0.011 \,\theta, \quad \theta \to 0.
\end{eqnarray*}
The local exact slope of the sequence $I_n^{(2)}$ as $\theta \to 0$ admits the representation
$$c_I^{(2)}(\theta)=(b^{(2)}_I(\theta))^2/(9\Delta_2^2) \sim 0.037\theta^2.$$

\noindent From \eqref{kul} the Kullback-Leibler "distance" for Makeham distribution satisfies
\begin{equation}
\label{KLM}
K_1(\theta) \sim \frac{\theta^2}{24}, \;\; \theta \to
0.
\end{equation}

\noindent Hence the local BE is
$$e^B(I^{(2)})=\lim_{\theta \to 0}\frac{c_I^{(2)}(\theta)}{2K_1(\theta)}= 0.448.$$

 The calculation
for other alternatives is quite similar, therefore we omit it and we present local Bahadur efficiencies in table \ref{fig: tab1}.

\begin{table}[!hhh]\centering
\bigskip
\caption{Local Bahadur efficiency for statistic $I_n^{(2)}$  }
\bigskip
\begin{tabular}{|c|c|}
\hline
Alternative & Efficiency\\
\hline
Makeham & 0.448\\
Weibull & 0.621\\
Gamma & 0.723\\
 EMNW(3) & 0.694\\
\hline
\end{tabular}
\label{fig: tab1}
\end{table}

\subsection{ Integral statistic $I_{n}^{(3)}$}

For $k=3$ from (\ref{psi_Ar}) and (\ref{sigma}) we get that the projection of the kernel  $\Psi_{3}(X,Y,Z,W)$ is equal to
\begin{gather}\label{psi_2}
\psi_3(s)=\frac{5}{8}e^{-s}-\frac{3}{4}e^{-2s}+\frac{1}{4}e^{-3s}-\frac{3}{32}e^{-s/2}-\frac{1}{12}e^{-s/3},
\end{gather}
\noindent and its variance is
\begin{gather*}
\Delta_3^2 = \int_{0}^{\infty} \psi_3^2 (s) e^{-s}ds =\frac{14591}{30750720}\approx
0.000474.
\end{gather*}

As in the previous case, according to Hoeffding's theorem, as $n \rightarrow \infty,$ the following convergence in distribution holds
$$\sqrt{n}I_n^{(3)} \stackrel{d}{\longrightarrow}{\cal
{N}}\left(0,\frac{14591}{1921920}\right).$$

 Regarding the large deviation asymptotics of the sequence $I_{n}^{(3)}$ under the null hypothesis,
we get exactly in the same manner as in the previous case:
\begin{theorem}
 For $a>0$ it holds
$$
\lim_{n\to \infty} n^{-1} \ln P_{H_0} ( I_n^{(3)}  >a) = - f (a),
$$
where the function $f$ is analytic for sufficiently small $a>0,$ moreover
\begin{gather}
f(a) \sim \frac{a^2}{32 \Delta_3^2}=  \frac{960960}{14591}a^2=65.86 a^2, \,
\, \mbox{as} \, \, a \to 0.\label{th1}
\end{gather}
\end{theorem}

 In this case the limit in probability under alternative $H_1$ is equal to
\begin{gather*}
b_I^{(3)}(\theta)=P_{\theta}(\max(X,Y, Z)<W)-P_{\theta}(X+\frac{Y}{2}+\frac{Z}{3}< W).
\end{gather*}
It is easy to show (\cite{NiPe}) that
$b_I^{(3)}(\theta) \sim 4\theta  \int_{0}^{\infty} \psi_3 (s)h(s)ds,
$
where again $h(x)=\frac{\partial}{\partial\theta}g_1(x,\theta)\mid _{\theta=0}$ and $\psi_3 (s)$ is the projection from (\ref{psi_2}).

For the Makeham alternative we have
\begin{eqnarray*}
b_I^{(3)}(\theta) &\sim& 4\theta
\int_{0}^{\infty}(\frac{5}{8}e^{-s}-\frac{3}{4}e^{-2s}+\frac{1}{4}e^{-3s}-\frac{3}{32}e^{-s/2}-\frac{1}{12}e^{-s/3})
e^{-s}(2-2e^{-s}-s)ds \\
&=&\frac{2}{105}\theta \approx
0.019 \,\theta, \quad \theta \to 0,
\end{eqnarray*}
and the local exact slope of the sequence $I_n^{(3)}$ as $\theta \to 0$ admits the representation
$$c_I^{(3)}(\theta)=(b^{(3)}_I(\theta))^2/(16\Delta_3^2) \sim 0.048\theta^2.$$

As previosly stated, the Kullback-Leibler "distance" satisfies the relation  (\ref{KLM}). Hence the local BE is equal to
$$e^B(I^{(3)})=\lim_{\theta \to 0}\frac{c_I^{(3)}(\theta)}{2K_1(\theta)}\approx 0.573.$$

We again omit the calculations
for other alternatives and we present local Bahadur efficiencies in table \ref{fig: tab2}.

\begin{table}[!hhh]\centering
\caption{Local Bahadur efficiency for $I_n^{(3)}$ }
\bigskip
\begin{tabular}{|c|c|}
\hline
Alternative & Efficiency\\
\hline
Makeham & 0.573\\
Weibull & 0.664\\
Gamma & 0.708\\
EMNW(3) & 0.799\\
\hline
\end{tabular}
\label{fig: tab2}
\end{table}

\bigskip

Using the MAPLE package we obtained maximal (with respect to $k$)  values of efficiencies against our four alternatives.
In table \ref{fig: tab3} we present the efficiencies from  tables \ref{fig: tab1} and \ref{fig: tab2} as well as the maximal values
we obtained.
\bigskip
\begin{table}[!hhh]\centering
\bigskip
\caption{Comparative table of local efficiencies for statistic $I_n^{(k)}$ }
\bigskip
\begin{tabular}{|c|c|c|c|}
\hline
Alternative & eff $k=2$ & eff $k=3$ & $\max_{k}$ eff  \\
\hline
Makeham & 0.448 & 0.573 & 0.875 for $k=14$\\
Weibull & 0.621 & 0.664 & 0.710 for $k=8$\\
Gamma & 0.723 & 0.708 & 0.723 for $k=2$\\
EMNW(3) &  0.694 & 0.799 & 0.885 for $k=6$ \\
\hline
\end{tabular}
\label{fig: tab3}
\end{table}

In table \ref{fig: tabsimint} we present the simulated powers for our four alternatives. The simulations have been performed for $n=100$ with 10000 replicates.




\begin{table}[htbp]
\centering
\bigskip
\caption{Simulated powers for statistic $I_n^{(k)}.$ }
\bigskip
\centering
\begin{tabular}{|c|c|c|c|c|c|}
  \hline
 Alternative & $\theta$ & $k$ & $\alpha=0.05$ & $\alpha=0.025$ & $\alpha=0.01$ \\
  \hline
   &0.5 & 2 & 0.1768& 0.1212& 0.0612  \\
   &0.5 & 3 & 0.2205&0.1306& 0.0706\\
  Makeham&0.5 & 4 & 0.2398&0.1532& 0.0772  \\
   &0.25 & 2 & 0.1091& 0.0653& 0.0294  \\
   &0.25 & 3 & 0.1171&0.0679& 0.0338\\
   &0.25 & 4 & 0.1392&0.0705& 0.0347  \\
   \hline
   &0.5 & 2 & 0.9963&0.9914& 0.9752  \\
   &0.5 & 3 & 0.9977&0.9942& 0.9839 \\
  Weibull&0.5 & 4 & 0.9987 &0.9965 &0.9864 \\
   &0.25 & 2 & 0.7166& 0.6456& 0.5049 \\
   &0.25 & 3 & 0.7626& 0.6456& 0.5049 \\
   &0.25 & 4 & 0.7940& 0.6813& 0.5309  \\
 \hline
   &0.5 & 2 & 0.8456&0.7736&0.6187 \\
   &0.5& 3 & 0.8453& 0.7528& 0.6198  \\
  Gamma&0.5 & 4 & 0.8528&0.7577& 0.6084  \\
   &0.25 & 2 & 0.4108& 0.3179& 0.1854  \\
   &0.25 & 3 & 0.4201&0.2940& 0.1836 \\
   &0.25 & 4 & 0.4323&0.3046& 0.1813  \\
 \hline
  &0.5&2& 0.9892 &0.9736& 0.9262\\
&0.5&3& 0.9841& 0.9591& 0.9097\\
 EMNW(3)&0.5&4&  0.9792& 0.9502& 0.8893\\
&0.25&2&  0.4476 &0.3454 &0.2098\\
&0.25&3&  0.4723& 0.3398& 0.2191\\
&0.25&4& 0.4820 &0.3577& 0.2173 \\
\hline

\end{tabular}
\label{fig: tabsimint}
\end{table}

\newpage
\section{ Kolmogorov-type statistic $D_n^{(k)}$}

In this section we consider the Kolmogorov-type statistic  (\ref{D_n_Ar}). For a fixed $t>0$ the expression $H_n^{(k)}(t)-G_n^{(k)}(t)$ is the $V$-statistic with the following kernel:

\begin{equation}
 \Xi_k(X_{1},X_{2},\ldots,X_{k};t)=\textbf{1}\{\max(X_{1},X_{2},\ldots,X_{k})<t\}-\frac{1}{k!}\sum\limits_{\pi(j_1, \ldots,j_k)}
 \textbf{1}\{\frac{X_{1}}{j_1}+\frac{X_{2}}{j_2}+\ldots+\frac{X_{k}}{j_k}<
 t\}.\nonumber
\end{equation}

Let $\xi_k (X_{1};t)$ be the projection of $\Xi_k(X_{1},X_{2},\ldots,X_{k};t)$ on $X_{1}$. Then
\begin{eqnarray}\label{ksikst}
\nonumber\xi_k (s;t)&=& E(\Xi_k(X_{1},X_{2},\ldots,X_{k};t)|X_{1}=s)\\
&=&\nonumber P\{ \max(s,X_{2},\ldots,X_{k})<t\}-\frac{1}{k!}\sum\limits_{\pi(j_1, \ldots,j_k)}
 P\{ \frac{s}{j_1}+\frac{X_{2}}{j_2}+\ldots+\frac{X_{k}}{j_k}<t\}\\
&=&\textbf{1} \{ s<t\}(F(t))^{k-1}-\frac{1}{k}\sum\limits_{j=1}^{k}\big[\textbf{1} \{ s<jt\} \big(1-\sum\limits_{\underset{i\neq j}{i=1}}^{k}\big(e^{-i(t-\frac{s}{j})}
\prod\limits_{\underset{h\neq i,j}{h=1}}^{k}\frac{h}{h-i}\big)\big)\big],
\end{eqnarray}
\noindent where $F(t)$ is d.f. of exponential distribution.
The calculation of variance for this projection in terms of $k$  is too complicated, therefore we calculate it only for particular cases.

\subsection{ Kolmogorov-type statistic $D_n^{(2)}$}

For $k=2$ from \eqref{ksikst} we get that the projection of the family of kernels $\Xi_2 (X,Y;t)$  is equal to
\begin{equation}\label{xi}
\xi_2 (s;t)=  \textbf{1}\{s<t\}F(t)-\frac12 \textbf{1}\{s<t\} F(2(t-s))-\frac12 \textbf{1}\{s< 2t\}F(t-s/2).
\end{equation}

Now we calculate the variances of these projections  $\delta_2^2(t)$ under $H_{0}.$ Elementary calculations show that
\begin{equation*}
\delta_2^2(t)=\frac{1}{3}e^{-t}-\frac{5}{4}e^{-2t}-\frac{1}{3}e^{-3t}
-\frac{1}{12}e^{-4t}-\frac23 e^{-3t/2}+2e^{-5t/2}+\frac12 t e^{-2t}.
\end{equation*}

\noindent Hence our family of kernels $\Xi_2 (X,Y;t)$  is non-degenerate as defined in \cite{Niki10} and besides
\begin{equation*}
\delta_2^2=\sup_{ t\geq 0} \delta_2^2(t)=0.02234.
\end{equation*}

\begin{figure}[h!]
\begin{center}
\includegraphics[scale=0.35]{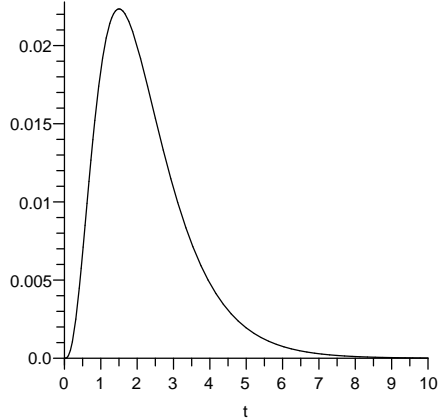}\caption{Plot of the function $\delta^2(t),$ $k=2.$ }
\end{center}
\end{figure}

Limiting distribution of the statistic $D_n^{(2)}$ is unknown. Using the methods of Silverman \cite{Silv}, one can show that the
$U$-empirical process
$$\eta_n^{(2)}(t) =\sqrt{n} \left(H_n^{(2)}(t) - G_n^{(2)}(t)\right), \ t\geq 0,
$$
weakly converges in $D(0,\infty)$ as $n \to \infty$ to certain centered Gaussian
process $\eta^{(2)}(t)$ with calculable covariance. Then the sequence of statistics
$\sqrt{n} D_n^{(2)}$ converges in distribution to the random variable   $\sup_{t\geq0} |\eta^{(2)}(t)|$ but it is impossible to find explicitly its distribution. Hence it is reasonable to determine the critical values for statistics  $D_n^{(2)}$ by simulation. Therefore in table \ref{fig: critic} we give the critical values for Kolmogorov-type statistic for $k=2$ and $k=3$ obtained via simulation.

\begin{table}[!hhh]\centering
\bigskip
\caption{Critical values for Kolmogorov type test $(n=100)$ }
\bigskip
\begin{tabular}{|c|c|c|c|c|}
\hline
$k$ & $\alpha=0.1$ & $\alpha=0.05$ & $\alpha=0.01$ &$\alpha=0.005$ \\
\hline
2 & 0.305 & 0.313 & 0.328 &0.334\\
3 & 0.446 & 0.455 & 0.473 &0.481\\
\hline
\end{tabular}
\label{fig: critic}
\end{table}

The family of kernels  $\{\Xi_2 (X,Y;t)\}, t\geq 0,$ is centered and bounded in the sense described in \cite{Niki10}. Applying the large deviation theorem for the supremum of the family of non-degenerate $U$- and $V$-statistics from \cite{Niki10} ,
we get the following result.
\begin{theorem}
For $a>0$ it holds
$$
\lim_{n\to \infty} n^{-1} \ln P_{H_0} ( D_n^{(2)} >a) = - f_2(a),
$$
where the function $f_2$ is continuous for sufficiently small $a>0,$ moreover
 $$
f_2(a) = (8 \delta_2^2)^{-1} a^2(1 + o(1)) \sim 5.595 a^2, \, \mbox{as}
\, \, a \to 0.
$$
\end{theorem}

\subsubsection{Local Bahadur efficiency of the statistic $D_n^{(2)}$}

According to Glivenko-Cantelli theorem
for $V$-statistics \cite{Jan} the limit in probability under the alternative for statistics $D_n^{(2)}$ is equal to
\begin{gather*}
b_D^{(2)}(\theta)= \sup_{t\geq 0} |b_D^{(2)}(t,\theta)|=\sup_{t\geq 0}
|P_{\theta}(\max(X,Y)<t)-P_{\theta}(X+\frac{Y}{2}< t)|.
\end{gather*}
Assuming the regularity of the alternative d.f., we can deduce
\begin{equation}\label{b_sup_2}
b_D^{(2)}(t,\theta) \sim 2\theta  \int_{0}^{\infty} \xi_2(s; t)h(s)ds,\hspace{3mm} \theta \to 0,
\end{equation}
 where again $h(x)=\frac{\partial}{\partial\theta}g(x,\theta)\mid _{\theta=0}$ and $\xi_2 (s;t)$ is the projection from (\ref{xi}).

We now proceed with calculation of local Bahadur efficiencies for our four alternatives.

For Makeham alternative from \eqref{b_sup_2} we get that
\begin{eqnarray*}
b_D^{(2)}(t,\theta) &\sim& \theta \bigg(2 \int_{0}^{t}F(t)e^{-s}(2-2e^{-s}-s)ds- \int_{0}^{t}F(2(t-s))e^{-s}(2-2e^{-s}-s)ds\\
 &&-\int_{0}^{2t}F(t-s/2)e^{-s}(2-2e^{-s}-s)ds\bigg)\\
 &=&\theta \big(\frac23 e^{-t}+(1-2t)e^{-2t}-2e^{-3t}+\frac13 e^{-4t}\big),\hspace{2mm} \theta
\to 0.
\end{eqnarray*}
Thus we have that
\begin{equation*}
\sup_{t>0}b_D^{(2)}(t,\theta)=
b_D^{(2)}(1.908,\theta) \sim 0.03055\; \theta, \, \theta \to 0.
\end{equation*}
\begin{figure}[h!]
\begin{center}
\includegraphics[scale=0.35]{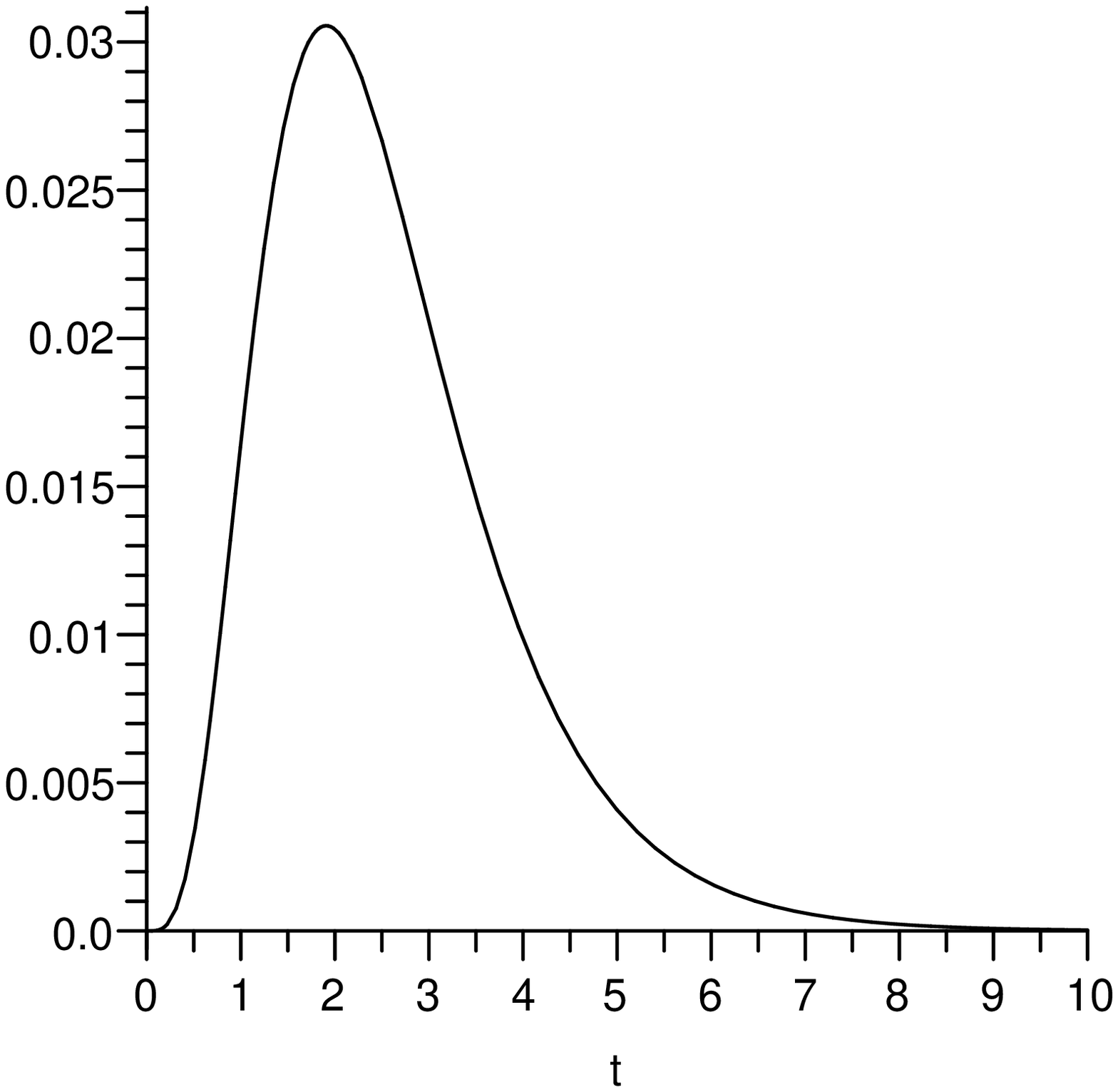}\caption{Plot of the function  $b_D^{(2)}(t,\theta), \mbox{ Makeham alternative}$ }
\end{center}
\end{figure}

The local exact slope of the sequence $D_n^{(2)}$ as $\theta \to 0$ satisfies
$$c_D^{(2)}(\theta)=(b^{(2)}_D(\theta))^2/(4 \delta_2^2)\sim 0.0104\;\theta^2.$$

Using $K_1(\theta)$ from (\ref{KLM}), we get that the local BE is equal to
$$e^B(D^{(2)})=\lim_{\theta \to 0}\frac{c_D^{(2)}(\theta)}{2K_1(\theta)}\approx 0.125.$$

For other alternatives the calculations are similar. Therefore we omit them and present their local Bahadur efficiencies in table \ref{fig: tab21}.

\begin{table}[!hhh]\centering
\caption{Local Bahadur efficiency for the statistic $D_n^{(2)}.$}
\bigskip
\begin{tabular}{|c|c|}
\hline
Alternative & Efficiency\\
\hline
Makeham & 0.125 \\
Weibull & 0.092 \\
Gamma & 0.093\\
EMNW(3)& 0.149\\
\hline
\end{tabular}
\label{fig: tab21}
\end{table}

We see that the efficiencies are very low, considerably lower than in case of other tests of exponentiality based on characterizations with the exception, apparently, of \cite{Niknik}. Probably this is related to intrinsic properties of Arnold-Villasenor characterization.

\subsection { Kolmogorov-type statistic $D_n^{(3)}$}

For $k=3$ from \eqref{ksikst} we get that the projection of the family of kernels $\Xi_3 (X,Y,Z;t)$  is equal to
\begin{eqnarray}\label{xi_3}
\nonumber \xi_3 (s;t)&=&  \textbf{1}\{x<t\}\big[ F^2(t)- F(2(t-x))+\frac23 F(3(t-x)) \big]- \textbf{1}\{x< 2t\}\big[\frac12F(t-x/2)\\
 &-&\frac16F(3(t-x/2))\big]- \textbf{1}\{x<3t\}\big[ \frac23 F(t-x/3)-\frac13 F(2(t-x/3))\big].
\end{eqnarray}

Now we calculate the variances of these projections  $\delta_3^2(t)$ under $H_{0}.$ We get that

\begin{eqnarray*}
\delta_3^2(t)&=&\frac{8}{15}e^{-t}+(\frac12 t-\frac{1}{24})e^{-2t}+(\frac{41}{9}-\frac43 t)e^{-3t}
-\frac{179}{210}e^{-4t}+\frac{113}{210}e^{-5t}-\frac{419}{2520}e^{-6t}\\
&-&\frac{14}{15} e^{-3t/2}+\frac{122}{35} e^{-5t/2}
-\frac{2}{3}e^{-7t/2}-\frac{2}{3}e^{-9t/2}-\frac{5}{7}e^{-5t/3}-\frac{5}{2}e^{-7t/3}+\frac{10}{7}e^{-8t/3}\\
&-& 4e^{-10t/3}
-2e^{-11t/3} +2e^{-13t/3}.
\end{eqnarray*}
The plot of this function is given in figure 3.

\begin{figure}[h!]
\begin{center}
\includegraphics[scale=0.35]{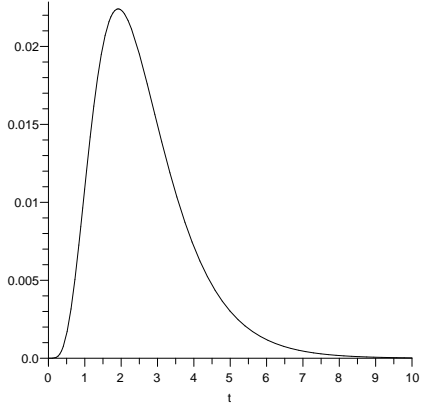}\caption{Plot of the function  $\delta_3^2(t).$  }
\end{center}
\label{fig: tabd3}
\end{figure}

\noindent Hence our family of kernels $\Xi_3 (X,Y,Z;t)$ is non-degenerate in the sense described in \cite{Niki10} and
\begin{equation*}
\delta_3^2=\sup_{ t\geq 0} \delta_3^2(t)=0.02241.
\end{equation*}

Using the same reasoning as in the case $D_n^{(2)}$ we conclude that it is impossible to find explicitly the limiting distribution of the statistic $D_n^{(3)}$.
The family of kernels  $\{\Xi_3 (X,Y,Z;t)\},$ $t\geq 0,$ is centered and bounded in the sense given in \cite{Niki10}. Applying the large deviation theorem for the supremum of the family of non-degenerate $U$- and $V$-statistics from \cite{Niki10},
we get the following result.
\begin{theorem}
For $a>0$ it holds
$$
\lim_{n\to \infty} n^{-1} \ln P_{H_0} ( D_n^{(3)} >a) = - f_3(a),
$$
where the function $f$ is continuous for sufficiently small $a>0,$ moreover
 $$
f_3(a) = (18 \delta_3^2)^{-1} a^2(1 + o(1)) \sim 2.479 a^2, \, \mbox{as}
\, \, a \to 0.
$$
\end{theorem}

\subsubsection{Local Bahadur efficiency of the statistic $D_n^{(3)}$}

In this case the limit in probability under the alternative, according to Glivenko-Cantelli theorem for $V$-statistics \cite{Jan}, is equal to
\begin{gather*}
b_D^{(3)}(\theta)= \sup_{t\geq 0}|b_D^{(3)}(t,\theta)|=\sup_{t\geq 0}
|P_{\theta}(\max(X,Y, Z)<t)-P_{\theta}(X+\frac{Y}{2}+\frac{Z}{3}< t)|.
\end{gather*}
It is not difficult to show that $b_D(t, \theta)$ for regular alternatives satisfies the relation
\begin{gather}\label{bd3}
b_D^{(3)}(t,\theta) \sim 3\theta  \int_{0}^{\infty} \xi_3 (s; t)h(s)ds,
\end{gather}
where $h(x)=\frac{\partial}{\partial\theta}g(x,\theta)\mid _{\theta=0}$, and $\xi_3 (s;t)$ is the projection from (\ref{xi_3}).

As in the previous sections we first calculate local BE for Makeham alternative. From \eqref{bd3} we get that
\begin{eqnarray*}
b_D^{(3)}(t,\theta) &\sim& \theta \bigg( \int_{0}^{t}\big[F^2(t)- F(2(t-s))+\frac23 F(3(t-s)) \big]e^{-s}(2-2e^{-s}-s)ds \\
&-&\int_{0}^{2t}\big[\frac12F(t-s/2) -\frac16F(3(t-s/2)) \big]e^{-s}(2-2e^{-s}-s)ds\\
&-&\int_{0}^{3t}\big[ \frac23 F(t-s/3)-\frac13 F(2(t-s/3)) \big]e^{-s}(2-2e^{-s}-s)ds\bigg)\\
&=&\theta\left(\frac85 e^{-t}+\big(\frac92-6t\big) e^{-2t}-8 e^{-3t}+2 e^{-4t}-\frac{1}{10} e^{-6t}\right), \ \theta
\to 0.
\end{eqnarray*}
Therefore we get that
\begin{equation*}
\sup_{t>0}b_D^{(3)}(t,\theta)=
b_D^{(3)}(2.087,\theta) \sim 0.0602\; \theta.
\end{equation*}

\begin{figure}[h!]
\begin{center}
\includegraphics[scale=0.35]{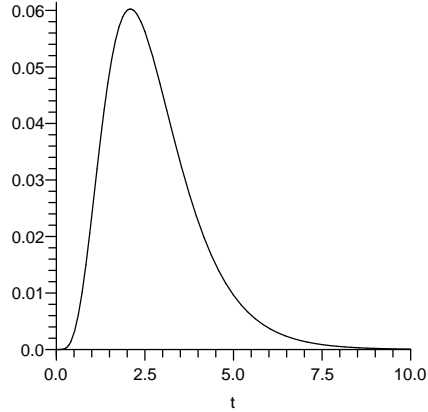}\caption{Plot of the function  $b_D^{(3)}(t,\theta), \mbox{ Makeham alteranative}$ }
\end{center}
\end{figure}

\noindent The local exact slope of the sequence $D_{n}^{(3)}$ as $\theta \to 0$
satisfies
\begin{equation}
c_D^{(3)}(\theta)=(b^{(3)}_D(\theta))^2/(9 \delta_3^2)\sim 0.018\;\theta^2,
 \end{equation}
and the local BE is equal to
$e^B(D^{(3)})=0.216.$ Omitting again the detailed calculations, we present in table \ref{fig: td3} the values of local Bahadur efficiency for our alternatives.

\bigskip

\begin{table}[!hhh]\centering
\bigskip
\caption{Local Bahadur efficiency for statistic $D_n^{(3)}$}
\bigskip
\begin{tabular}{|c|c|}
\hline
Alternative  & Efficiency\\
\hline
Makeham & 0.216 \\
Weibull & 0.152 \\
Gamma & 0.138\\
EMNW(3) & 0.230\\
\hline
\end{tabular}
\label{fig: td3}
\end{table}

We see that these efficiencies are slightly better than in the previous case, but still rather low. In table \ref{fig: tabsimks} we present the simulated powers for our four alternatives.
Again the simulations have been performed for $n=100$ with 10000 replicates.





\begin{table}[htbp]
\centering
\bigskip
\caption{Simulated powers for statistic $D_n^{(k)}$ }
\bigskip
\centering
\begin{tabular}{|c|c|c|c|c|c|}
  \hline
  Alternative & $\theta$ & $k$ & $\alpha=0.05$ & $\alpha=0.025$ & $\alpha=0.01$  \\
  \hline
   &0.5 & 2 &  0.0885& 0.0472 &0.0221  \\
   &0.5 & 3 & 0.1027& 0.0609 &0.0246\\
  Makeham&0.5 & 4 & 0.1136& 0.0681& 0.0304  \\
   &0.25 & 2 & 0.0669& 0.0315& 0.0154  \\
   &0.25 & 3 & 0.0724& 0.0399& 0.0164 \\
   &0.25 & 4 & 0.0842& 0.0475& 0.0206 \\
  \hline
   &0.5 & 2&0.6967& 0.5721& 0.4423 \\
   &0.5 & 3 & 0.8194& 0.7431& 0.6006  \\
  Weibull&0.5 & 4 & 0.8903& 0.8287& 0.7190\\
   &0.25 & 2 & 0.2969 &0.1964& 0.1169 \\
   &0.25 & 3 & 0.3698 & 0.2745& 0.1566 \\
   &0.25 & 4 & 0.4286 & 0.3308& 0.2054   \\
  \hline
   &0.5 & 2&  0.4146& 0.2901& 0.1849 \\
   &0.5& 3&0.5026& 0.3887& 0.2405   \\
  Gamma&0.5 & 4& 0.5555& 0.4433& 0.3006 \\
   &0.25 & 2& 0.1852& 0.1135& 0.0630   \\
   &0.25 & 3& 0.2163& 0.1437& 0.0695  \\
   &0.25 & 4&0.2406& 0.1628& 0.0841   \\
  \hline
 &0.5&2&0.7083 &0.5769& 0.4352\\
&0.5&3&0.7918& 0.6936& 0.5294\\
EMNW(3)&0.5&4& 0.8409& 0.7581& 0.6121\\
 &0.25&2&0.2080& 0.1294& 0.0718\\
 &0.25&3& 0.2456& 0.1658 &0.0817\\
 &0.25&4&  0.2849& 0.1964& 0.1083\\
   \hline

\end{tabular}
\label{fig: tabsimks}
\end{table}

\newpage
\section{Conditions of local asymptotic optimality} \label{localop}

The efficiencies of our tests for standard alternatives are far from maximal ones. Nevertheless, there exist special alternatives (we call them {\it most favorable})
for which our sequences of statistics $I_n^{(k)}$ and $D_n^{(k)}$ are locally asymptotically optimal (LAO) in Bahadur sense (see general theory in \cite[Ch.6]{Nik}).
In this section we describe the local structure of such
alternatives, for which the given statistic has maximal possible local efficiency, so that the relation
$$
c_T(\theta) \sim 2 K(\theta),\, \theta \to 0,
$$
holds (see \cite{Bahadur}, \cite{Nik}, \cite{NikTchir}, \cite{NiPe}).
Such alternatives form the so-called domain of LAO for the given sequence of statistics $\{T_n\}$.

Denote by $\cal G$
the class of densities $ g(\cdot \ ,\theta)$ with the d.f.'s $G(\cdot \ ,\theta)$. Define the functions
\begin{gather*}
H(x)=\frac{\partial}{\partial\theta}G(x,\theta)\mid
_{\theta=0},\quad
h(x)=\frac{\partial}{\partial\theta}g(x,\theta)\mid _{\theta=0}.
\end{gather*}

\noindent Suppose also that for $G$ from $\cal G$ the following regularity conditions hold:
\begin{gather*}
h(x)=H'(x), \,  x \geq 0, \quad \int_0^\infty h^2(x)e^{x}dx <  \infty  , \\
\frac{\partial}{\partial\theta}\int_0^\infty x g(x,\theta)dx \mid
_{\theta=0} \ = \ \int_0^\infty x h(x)dx.
\end{gather*}

It is easy to show, see also \cite{NikTchir}, that under these conditions
$$ 2K(\theta)\sim \bigg[\int_0^\infty h^2(x)e^{x}dx -\bigg(\int_0^\infty x h(x)dx\bigg)^2\bigg] \theta^2,\, \theta \to 0.
$$

It can be shown that for the statistic (\ref{I_n_Ar}) holds
 \begin{equation*}
 b_I^{(k)}(\theta) \sim (k+1)\theta  \int_{0}^{\infty} \psi_k (x)h(x)dx, \, \theta \to 0.
 \end{equation*}
  Let us introduce the auxiliary function
\begin{equation}
\label{h0}
h_0(x) = h(x) - (x-1)\exp(-x)\int_0^\infty u h(u) du.
\end{equation}
It is straightforward that

\begin{eqnarray}\label{LaoD}
\int_0^\infty h^2(x)e^{x}dx -\bigg(\int_0^\infty x h(x)dx\bigg)^2 &=& \int_0^\infty h_0^2(x) e^{x} dx,\\
\nonumber \int_{0}^{\infty} \psi_k (x)h(x)dx &=& \int_{0}^{\infty} \psi_k (x)h_0(x)dx.
\end{eqnarray}

Consequently the local BE takes the form
\begin{eqnarray*}
e^B(I_{n}^{(k)})&=& \lim_{\theta \to 0} \frac{(b_I^{(k)}(\theta))^2}{ 2(k+1)^2 \Delta_k^2 K(\theta)}\\
&=& \bigg( \int_{0}^{\infty} \psi_k (x)h_0(x)dx\bigg)^2/\bigg(
\int_{0}^{\infty}\psi_k^2(x) e^{-x}dx \cdot  \int_0^\infty h_0^2(x)e^{x}dx
 \bigg).
\end{eqnarray*}

The local Bahadur asymptotic optimality means that the expression
on the right-hand side is equal to 1. It follows from Cauchy-Schwarz inequality
(see also \cite{NiPe})
that this is satisfied if
$h_0(x)=C_1 e^{-x}\psi(x)$ for some constant $C_1>0,$
so that $h(x) = e^{-x}(C_1\psi(x)+
C_2(x-1))$ for some constants $C_1>0$ and $C_2.$
Such distributions constitute the LAO domain in the class $\cal G$.

The simplest examples of
such alternative densities  $g(x,\theta)$  for small $\theta > 0$ are given in table \ref{fig: lao1}.

\bigskip

\begin{table}[!hhh]\centering
\caption{Most favorable alternatives for $I_n^{(k)}.$}
\bigskip
\begin{tabular}{|c|l|}
\hline
& \\
 & Alternative density  $g(x, \theta)$ as $\theta \to +0, \, x \geq 0$\\
\hline
& \\
$k=2$ &  $g(x,\theta)=e^{-x}\big(1+\frac{\theta}{3}
(\frac{4}{3}e^{-x}-e^{-2x}-\frac{1}{2}e^{-x/2})\big)$\\
\hline
& \\
$k=3$& $g(x,\theta)=e^{-x}\big(1+\frac{\theta}{4}
(\frac{5}{2}e^{-x}-3e^{-2x}+e^{-3x}-\frac{3}{8}e^{-x/2}-\frac{1}{3}e^{-x/3})\big)$\\
\hline
\end{tabular}
\label{fig: lao1}
\end{table}

Let us now consider the Kolmogorov-type statistic (\ref{D_n_Ar}). It can be shown that
 \begin{equation*}
 b_D^{(k)}(\theta) \sim k\theta  \int_{0}^{\infty} \xi_k (x;t)h(x)dx, \, \theta \to 0.
 \end{equation*}

\noindent For $h_{0}(x)$ defined in \eqref{h0}, besides \eqref{LaoD}, also holds

\begin{equation*}
 \int_{0}^{\infty} \xi_k (x;t)h(x)dx = \int_{0}^{\infty} \xi_k (x;t)h_0(x)dx.
\end{equation*}

\noindent In this case the efficiency is equal to
\begin{eqnarray*}
e^B (D_{n}^{(k)})&=& \lim_{\theta \to 0}\frac{ (b^{(k)}_{D}(\theta))^{2}} {\sup_{t\geq 0}\left(2k^2
\delta_k^2(t)\right) K(\theta) }\\
&=& \sup_{t\geq 0}\left( \int_{0}^{\infty}
 \xi_k (x;t)h_0(x)dx\right)^2 / \ \sup_{t\geq 0} \left(
\int_{0}^{\infty}\xi_k^2 (x;t) e^{-x}dx \cdot \int_0^\infty h_0^2 e^{x} dx\right).
\end{eqnarray*}

 From Cauchy-Schwarz inequality we obtain that efficiency is equal to 1 if
$ h(x)=e^{-x}\big(C_1\xi_k(x; t_0)+ C_2(x-1)\big)$
for $t_0= {\rm argmax}_{t\geq0}
\delta_k^2(t) $ and some constants $C_1>0$ and $C_2.$
The alternative densities having such function $h(x)$ form the domain of LAO in the
corresponding class.

The simplest examples are given in table \ref{fig: tlao2}. To facilitate the presentation, we denote:
\begin{eqnarray*}
t_0&=& \underset{t\geq0}{\rm argmax} \bigg(
\frac{1}{3}e^{-t}-\frac{5}{4}e^{-2t}-\frac{1}{3}e^{-3t}
-\frac{1}{12}e^{-4t}-\frac23 e^{-3t/2}+2e^{-5t/2}+\frac12 t e^{-2t}\bigg)\approx 1.502;\\
t_1&=& \underset{t\geq0} {\rm argmax} \bigg[ \frac{8}{15}e^{-t}+(\frac12 t-\frac{1}{24})e^{-2t}+(\frac{41}{9}-\frac43 t)e^{-3t}
-\frac{179}{210}e^{-4t}+ \frac{113}{210}e^{-5t}\\
&-&\frac{419}{2520}e^{-6t}-\frac{14}{15} e^{-3t/2}+\frac{122}{35} e^{-5t/2}-
\frac{2}{3}e^{-7t/2}-\frac{2}{3}e^{-9t/2}- \frac{5}{7}e^{-5t/3}-\frac{5}{2}e^{-7t/3}\\
&+&\frac{10}{7}e^{-8t/3}
-4e^{-10t/3}-2e^{-11t/3} +2e^{-13t/3} \bigg]\approx 1.919.
\end{eqnarray*}

\medskip

\begin{table}[!hhh]\centering
\caption{Most favorable alternatives for $D_n^{(k)}$}
\bigskip
\begin{tabular}{|c|l|}
\hline
& \\
 & Alternative densities $g(x, \theta)$ as $\theta \to +0, \, x \geq 0$\\
\hline
& \\
$k=2$ & $g(x,\theta)=e^{-x}\bigg(1+\theta \cdot
\textbf{1}\{x<t_0\}(1-e^{-t_0})-$\\
 & \qquad \qquad $-
\frac12 \theta \cdot \big(\textbf{1}\{x<t_0\}(1-e^{-2(t_0-x)}) + \textbf{1}\{x< 2t_0\}(1-e^{-(t_0-x/2)}) \big)\bigg)$\\
\hline
& \\
$k=3$ & $g(x,\theta)=e^{-x}\bigg(1+\theta \cdot
\textbf{1}\{x<t_1\}\big[(1-e^{-t_1})^2 +e^{-2(t_1-x)} -\frac23 e^{-3(t_1-x)}-\frac13 \big]$\\
 & \qquad \qquad \qquad  $-\frac13 \theta \cdot\textbf{1}\{x< 2t_1\}\big[1-\frac32e^{-(t_1-x/2)}+\frac12e^{-3(t_1-x/2)}\big] $\\
 & \qquad \qquad \qquad$-\frac13 \theta \cdot \textbf{1}\{x<3t_1\}\big[ 1-2e^{-(t_1-x/3)}+e^{-2(t_1-x/3)}    \big]\bigg)$\\
\hline
\end{tabular}
\label{fig: tlao2}
\end{table}

\section{Discussion}

In this paper we have proposed two families of asymptotic tests of exponentiality based on recent characterization of exponentiality by Arnold and Villasenor \cite{Ar}. The integral test statistics $I_n^{(k)}$ are asymptotically normal and have reasonably simple form which can be easily computed for small $k.$ They are consistent for many common alternatives and have local Bahadur efficiency around 0.5  - 0.7. There exist also special (most favorable) alternatives described in the section \ref{localop} for which the integral statistics are locally asymptotically optimal in this sense.

We also obtained via simulation the power of new integral statistics for chosen alternatives. In theory, the ordering
of tests by power is linked more closely to Hodges-Lehmann efficiency \cite{Nik}, and should not coincide with the ordering by local Bahadur efficiency. Nevertheless, we observe tolerable correspondence of test quality according to both criteria with the exception of Gamma and Weibull distribution. In whole we can recommend new integral tests of exponentiality as additional and auxiliary tests of exponentiality, especially when one is trying to reject exponentiality in a specific example using a "battery" of statistical tests.

In the case of Kolmogorov type tests the values of local Bahadur efficiency turned  out to be rather low for common alternatives, and the simulated powers (which are slightly more optimistic) do not change somewhat disadvantageous regard to new tests of exponentiality of supremum type. Probably it is closely related to intrinsic properties of Arnold-Villasenor characterization. However, even these tests, in virtue of their consistency, can be of some use in statistical research, especially when the (unknown) alternative is close to the most favorable one.

\section{Acknowledgement}

The authors express their deep gratitude to Prof. George Yanev who kindly sent them the files of his papers \cite{Chakra} and \cite{Yanev}.


\end{document}